\documentclass{article}
\usepackage{amsmath,amsthm,amsopn,amstext,amscd,amsfonts,amssymb}
\usepackage{natbib}

\def\diag{\mathop{\rm diag}\nolimits}
\def\tr{\mathop{\rm tr}\nolimits}

\def\build#1#2#3{\mathrel{\mathop{#1}\limits^{#2}_{#3}}}
\def\sign{\mathop{\rm sign}\nolimits}
\def\vec{\mathop{\rm vec}\nolimits}

\def\etr{\mathop{\rm etr}\nolimits}
\newcommand {\findemo}{\hfill \, $\Box$ \,\\[2ex]}

\renewenvironment{abstract}
                 {\vspace{6pt}
                  \begin{center}
                  \begin{minipage}{5in}
                  \centerline{\textbf{Abstract}}
                  \noindent\ignorespaces
                 }
                 {\end{minipage}\end{center}}
\newtheorem{thm}{\textbf{Theorem}}[section]
\newtheorem{cor}{\textbf{Corollary}}[section]
\newtheorem{lem}{\textbf{Lemma}}[section]
\theoremstyle{definition}

\title{\Large \textbf{Shape theory via SVD decomposition I}}
\author{
  \textbf{Jos\'e A. D\'{\i}az-Garc\'{\i}a} \thanks{Corresponding author\newline
    }\\
  {\normalsize Department of Statistics and Computation} \\
  {\normalsize Universidad Aut\'onoma Agraria Antonio Narro}\\
  {\normalsize 25350 Buenavista, Saltillo, Coahuila, Mexico} \\
  {\normalsize E-mail: jadiaz@uaaan.mx} \\[2ex]
  \textbf{Francisco J. Caro-Lopera} \\
  {\normalsize Department of Basic Sciences} \\
  {\normalsize Universidad de Medell\'{\i}n} \\
  {\normalsize Carrera 87 No.30-65, of. 5-103}\\
  {\normalsize Medell\'{\i}n, Colombia}\\
  {\normalsize E-mail: fjcaro@udem.edu.co}\\
}
\date{}
\begin{document}
\maketitle

\begin{abstract}
This work finds the non isotropic noncentral elliptical shape distributions via SVD
decomposition in the context of zonal polynomials, avoiding the invariant polynomials
and the open problems for their computation. The new shape distributions are easily
computable and then the inference procedure is based on exact densities instead of
the published approximations and asymptotic densities of isotropic models. An
application of the technique is illustrated with a classical landmark data in
Biology, for this, three models are proposed, the usual Gaussian and two non
Gaussian; the best one is chosen by using  a modified BIC criterion.
\end{abstract}

\section{Introduction}\label{sec:Intro}

The multivariate statistical theory of shape has been studied deeply in the last two
decades (\citet{GM93}, \citet{dgm:97}, \citet{dgr:03}, \citet{DM98} and the
references there in, \citet{Caro2009}, among many others.  Most of the works are
supported by important restrictions (isotropy) for the covariance matrix and models
(Normal) in order to obtain known polynomials (zonal polynomials). A sort of
approaches are given for shape theory, via QR (\citet{GM93}), SVD (\citet{LK93},
\citet{g:91}, \citet{dgr:03}), affine (\citet{GM93}, \citet{dgr:03},
\citet{Caro2009}).

Avoiding the restrictions of isotropy and normality carry some problems, because
integration over Euclidean or affine transformations lead to the apparition of
invariant polynomials of \citet{d:80} which can not computed for large degrees.

This work finds a sequence of transformations which let the construction of shape
densities via the singular value decomposition and based on a non restricted non
central and non isotropy elliptical model. The resulting densities avoids the
invariant polynomials and they are set in terms of series of zonal polynomials which
can be computed by suitable modifications of the existing algorithms for
hypergeometric series (\citet{KE06}).

The work is structured as follows, the main principle and the size and shape
distribution is given in section \ref{sec:SVDsizeandshape}, then the shape density is
obtained in section \ref{sec:svdreflectionshape}; the associated excluding reflection
densities are considered in section \ref{sec:SVDexcludingreflection} and finally some
particular models are derived as corollaries in section \ref{sec:particularmodels}
which also presents an application in mouse vertebra by studying three models (the
usual Gaussian and two non Gaussian) with the modified BIC criterion.

\section{Main principle and SVD size-and-shape density}\label{sec:SVDsizeandshape}

It is known that the shape of an object is all geometrical information that remains
after filtering out translation, rotation and scale information of an original figure
(represented by a matrix $\mathbf{X}$) comprised in $N$ landmarks in $K$ dimensions.
So, we say that two figures, $\mathbf{X}_{1}:N\times K$ and $\mathbf{X}_{2}:N\times
K$ have the same shape if they are related by a special  similarity transformation
$\mathbf{X}_{2}=\beta \mathbf{X}_{1}\mathbf{H}+\mathbf{1}_{N}\boldsymbol{\gamma}'$,
where $\mathbf{H}:K\times K\in SO(K)$ (the rotation), $\boldsymbol{\gamma}:K\times 1$
(the translation), $\mathbf{1}_{N}:N\times 1,$ $\mathbf{1}_{N}=(1,1,\ldots,1)'$, and
$\beta>0$ (the scale). Thus, in this context, the shape of a matrix $\mathbf{X}$ is
all the geometrical information about $\mathbf{X}$ that is invariant under Euclidean
similarity transformations.

Now, multivariate statistical theory of shape compares shapes of objects in presence
of  randomness, so if we assume that  a figure $\mathbf{X}$, comprised in $N$
landmarks in $K$, follows an elliptical distribution $\mathbf{X} \sim \mathcal{E}_{N
\times K} (\boldsymbol{\mu}_{{}_{\mathbf{X}}},
\boldsymbol{\Sigma}_{{}_{\mathbf{X}}},\boldsymbol{\Theta}, h)$, it is of interest to
remove translation, scaling, rotation from $\mathbf{X}$. Clearly, the sequence
$\mathbf{L}\mathbf{X}=\mathbf{Y}=\mathbf{H}'\mathbf{D}\mathbf{P}=r\mathbf{W}(\mathbf{u})\mathbf{P}$
removes the translation (by a sub Helmert matrix $\mathbf{L}$, for example), the
rotation (by the SVD of $\mathbf{Y}$) and the scale (by dividing for the norm of
$\mathbf{Y}$). In order to obtain the density of $\mathbf{W}$ we need to integrate
over the similarity group;  it is easy to see that the elliptical assumption lead to
the product of two traces which irremediably expands in terms of invariant
polynomials of two matrix arguments (\citet{d:80}), and the shape densities are not
computable for large degrees.

So the classical statistical multivariate analysis  restricts the models for the
original landmark data in order to obtain densities which are expanded in terms of
studied polynomials such as the zonal polynomials which are computable
(\citet{GM93}),  otherwise, as we proved in the last sentence, the densities involve
non computable polynomials for large degrees.

From the practical point of view the restrictions affect  the applications; i.e., the
isotropic assumption $\boldsymbol{\Theta}=\mathbf{I}_{K}$ for an elliptical shape
model of the form
$$
   \mathbf{X} \sim \mathcal{E}_{N \times K} (\boldsymbol{\mu}_{{}_{\mathbf{X}}},
        \boldsymbol{\Sigma}_{{}_{\mathbf{X}}},\boldsymbol{\Theta}, h),
$$
restricts substantially the correlations of the landmarks in the figure,  specially
in objects with symmetries as in the case of mouse vertebra for example, among many
others (\citet{DM98}).  So, we expect the non isotropic model, with any positive
definite matrix $\boldsymbol{\Theta}$, as the best model for considering all the
possible correlations among the anatomical (geometrical o mathematical) points.
However, using the classical approach of the published literature of shape (see for
example \citet{GM93}) under the non isotropic model, we obtain immediately invariant
polynomials, which can not be computed at this time for large degrees.

In order to avoid this problem  consider the following procedure: Let
$$
   \mathbf{X} \sim \mathcal{E}_{N \times K} (\boldsymbol{\mu}_{{}_{\mathbf{X}}},
   \boldsymbol{\Sigma}_{{}_{\mathbf{X}}},\boldsymbol{\Theta}, h),
$$
if $\boldsymbol{\Theta}^{1/2}$ is the positive definite square root of the matrix
$\boldsymbol{\Theta}$, i .e. $\boldsymbol{\Theta} = (\boldsymbol{\Theta}^{1/2})^{2}$,
with $\boldsymbol{\Theta}^{1/2}:$ $K \times K$, \citet[p. 11]{gv:93}, and noting that
$$
  \mathbf{X} \boldsymbol{\Theta}^{-1} \mathbf{X}' = \mathbf{X} (\boldsymbol{\Theta}^{-1/2}
  \boldsymbol{\Theta}^{-1/2})^{-1}\mathbf{X}' = \mathbf{X} \boldsymbol{\Theta}^{-1/2}
  (\mathbf{X} \boldsymbol{\Theta}^{-1/2})' = \mathbf{Z}\mathbf{Z}',
$$ where
$$
\mathbf{Z} = \mathbf{X} \boldsymbol{\Theta}^{-1/2},
$$
then
$$
  \mathbf{Z} \sim \mathcal{E}_{N \times K}(\boldsymbol{\mu}_{{}_{\mathbf{Z}}},
  \boldsymbol{\Sigma}_{{}_{\mathbf{X}}}, \mathbf{I}_{K}, h),
$$
with $\boldsymbol{\mu}_{{}_{\mathbf{Z}}} = \boldsymbol{\mu}_{{}_{\mathbf{X}}}
\boldsymbol{\Theta}^{-1/2}$, (see \citet[p. 20]{gv:93}).

And we arrive at the classical starting point in shape theory where the original
landmark matrix is replaced by $\mathbf{Z} = \mathbf{X} \boldsymbol{\Theta}^{-1/2}$.
Then we can proceed as usual, removing from $\mathbf{Z}$, translation, scale,
rotation and/or reflection in order to obtain the shape of $\mathbf{Z}$ (or
$\mathbf{X}$) via the SVD decomposition, for example.

The SVD decomposition  has two version in shape theory, \citet{g:91} and
\citet{LK93}, we focus in this paper on Goodall's approach.

Let $n = \min(m,K)$, $\mathbf{Y} = \mathbf{H}'\mathbf{D}\mathbf{P}$ be the
nonsingular part of the SVD, where $\mathbf{H}:n \times m$, $\mathbf{H} \in
\mathcal{V}_{n,m}$ represents the Stiefel manifold, $\mathbf{D} = \diag (D_{1},
D_{2}, \ldots, D_{n})$ with $D_{1} \geq D_{2} \geq \cdots \geq D_{n}$ and
$\mathbf{P}:n \times K$, $\mathbf{P} \in \mathcal{V}_{n,K}$.

Thus the SVD shape coordinates ${\bf u}$ of $\mathbf{X}$ may be found by the
following procedure
\begin{equation}
\label{svd1}
  \mathbf{L}\mathbf{X} \boldsymbol{\Theta}^{-1/2}= \mathbf{L}\mathbf{Z} = \mathbf{Y} =
  \mathbf{H}'\mathbf{D}\mathbf{P} = r \mathbf{W}\mathbf{P} = r \mathbf{W}({\bf u})\mathbf{P}
\end{equation}
where the SVD shape coordinate system is given by $\mathbf{H}'\mathbf{D}$ (\citet
[pp. 296-298]{g:91}) and $r =||\mathbf{H}\mathbf{D}|| = (\tr \mathbf{D} \mathbf{H}
\mathbf{H}'\mathbf{D})^{1/2} = (\tr \mathbf{D}^{2})^{1/2} = ||\mathbf{D}|| =
||\mathbf{Y}||$. Before defining $\mathbf{W}$ and ${\bf u}$, note that when $n = K $
two cases may be distinguished.
\begin{enumerate}
   \item $\mathbf{P}$ includes reflection, $\mathbf{P} \in \mathcal{O}(k)$, $|\mathbf{P}| =
        \pm 1$, $\mathbf{D}_{K} \geq 0$ and $(\mathbf{H},\mathbf{D})$, written $(\mathbf{H},
        \mathbf{D})^{R}$ for definiteness, contains reflection SVD shape co-ordinates.
   \item $\mathbf{P}$ excludes reflection, $\mathbf{P} \in \mathcal{SO}(K)$, $|\mathbf{P}| =
        + 1$, $|D_{K}| \geq 0$, $\sign(D_{m}) = \sign|\mathbf{X}|$ and $(\mathbf{H},\mathbf{D})$,
        may be written $(\mathbf{H},\mathbf{D})^{NR}$ for definiteness.
\end{enumerate}
Now the SVD shape matrix $\mathbf{W}$ is obtained by dividing the
$\mathbf{H}'\mathbf{D}$ matrix by $r$, when  $\mathbf{W}$ may include or exclude
reflection, in which case we obtain, respectively, $\mathbf{W}^{R} =
(\mathbf{H}'\mathbf{D})^{R}/r$ or $\mathbf{W}^{NR} = (\mathbf{H}'\mathbf{D})^{NR}/r$.
Finally ${\bf u}$ is composed of the $mn-1$ generalized polar coordinates.

Our interest now lies in finding the corresponding densities associated with the
process described in (\ref{svd1}). Thus we obtain the joint density of
$(\mathbf{H},\mathbf{D})$ and the density of $\mathbf{W}({\bf u})$.

In order to obtain the size and shape density we need some integrals involving zonal
polynomials, extending \citet[eq. (22)]{JAT64}.

\begin{lem}\label{lem:geneq22James64}
Let $\mathbf{X}: K \times n$, $\mathbf{Y}: K \times K$ and $\mathbf{H} \in
\mathcal{V}_{n,K}$. Then
\begin{enumerate}
  \item
  $$
    \int_{\mathbf{H} \in \mathcal{V}_{n,K}} [\tr(\mathbf{Y} + \mathbf{X}\mathbf{H})]^{p} (\mathbf{H}d\mathbf{H}')
    = \displaystyle\frac{2^{n}
    \pi^{Kn/2}}{\Gamma_{n}[\frac{1}{2} K]} \sum_{f=0}^{\infty} \sum_{\lambda}
    \displaystyle\frac{(p)_{2f} (\tr \mathbf{Y})^{p-2f}}{(\frac{1}{2} K)_{\lambda}}
    \displaystyle\frac{C_{\lambda}(\frac{1}{4} \mathbf{X}\mathbf{X}')}{f!}
  $$
  where $|(\tr \mathbf{Y})^{-1} \tr \mathbf{X}\mathbf{H}|<1$ and $\tr \mathbf{Y} \neq 0$.
  \item
  $
    \displaystyle \int_{\mathbf{H} \in \mathcal{V}_{n,K}} \tr(\mathbf{Y} + \mathbf{X}\mathbf{H})
    \etr\{r(\mathbf{Y} + \mathbf{X}\mathbf{H})\}(\mathbf{H}d\mathbf{H}') =
  $
  \par \noindent \hfill
    \hbox{$\displaystyle\frac{2^{n}
    \pi^{Kn/2}}{\Gamma_{n}[\frac{1}{2} K]} \etr\{r \mathbf{Y}\} \left \{ \tr \mathbf{Y}
    {}_{0}F_{1} (\frac{1}{2} K; \displaystyle\frac {r^{2}}{4} \mathbf{X}\mathbf{X}') + \sum_{f=0}^{\infty}
    \sum_{\lambda} \displaystyle\frac{(f + \frac{1}{2})}{(\frac{1}{2} K)_{\lambda}}
    \displaystyle\frac{C_{\lambda}(\frac{1}{4} \mathbf{X}\mathbf{X}')}{f!} \right \},$}
  \par \noindent
\end{enumerate}
where $p \in \Re$, $r \in \Re$, $C_{\kappa}(\mathbf{B})$ are the
zonal polynomials of $\mathbf{B}$ corresponding to the partition
$\kappa=(f_{1},\ldots f_{p})$ of $f$, with $\sum_{i=1}^{p}f_{i}=f$; and
$(a)_{\kappa}=\prod_{i=1}(a-(j-1)/2)_{f_{j}}$, $(a)_{f}=a(a+1)\cdots (a+f-1)$, are
the generalized hypergeometric coefficients and ${}_{0}F_{1}$ is the Bessel function,
\citet{JAT64}.
\end{lem}

\textit{Proof.}
\begin{enumerate}
   \item From Lemma 9.5.3 \citet[Lemma 9.5.3, p. 397]{MR1982} we have
      $$
        \int_{\mathbf{H} \in \mathcal{V}_{n,K}} [\tr(\mathbf{Y} + \mathbf{X}\mathbf{H})]^{p}
        (\mathbf{H}d\mathbf{H}') = \frac{2^{n}
        \pi^{Kn/2}}{\Gamma_{n}[\frac{1}{2} K]} \int_{\mathcal{O}(K)}[\tr(\mathbf{Y} +
        \mathbf{X}\mathbf{H})]^{p} (d\mathbf{H}).
      $$
      Furthermore, for $\tr \mathbf{Y} \neq 0$ and $|(\tr \mathbf{Y})^{-1} \tr \mathbf{X}\mathbf{H}| <1$
      $$
        [\tr(\mathbf{Y} + \mathbf{X}\mathbf{H})]^{p} = (\tr \mathbf{Y})^{p} \sum_{f=0}^{\infty}
        \frac{(p)_{f}}{f!}(\tr \mathbf{Y})^{-f}(\tr \mathbf{X}\mathbf{H})^{f}.
      $$
      Now from \citet[eqs. (46) and (22)]{JAT64}) it follows that\\[2ex]
      $
       \displaystyle \int_{\mathbf{H} \in \mathcal{V}_{n,K}} [\tr(\mathbf{Y} +
       \mathbf{X}\mathbf{H})]^{p} (\mathbf{H}d\mathbf{H}') =
      $
      \par \noindent \hfill
      \hbox{$\displaystyle\frac{2^{n}
        \pi^{Kn/2}}{\Gamma_{n}[\frac{1}{2} K]} \sum_{f=0}^{\infty}
        \sum_{\lambda} \frac{(p)_{2f} (\tr \mathbf{Y})^{-2f}}{(2f)!}
        \frac{(\frac{1}{2})_{f}}{(\frac{1}{2} K)_{\lambda}} C_{\lambda}(\mathbf{X}\mathbf{X}'),$}
      \par \noindent\\
      the result follows, noting that $(\frac{1}{2})_{f}/(2f)! = 1/(4^{f} f!)$ and that $C_{\lambda}
      (a\mathbf{X}\mathbf{X}') = a^{f} C_{\lambda}(\mathbf{X}\mathbf{X}')$.
   \item This follows by expanding the exponentials in series of powers and by applying (22)
        and (27) from \citet{JAT64}.\findemo
\end{enumerate}

Now, the jacobian of the corresponding decomposition is provided next:

\begin{lem}\label{lem:SVDjacobian}
Let be $\mathbf{Y}:N-1\times K$, then there exist $\mathbf{V}\in
\mathcal{V}_{n,N-1}$, $\mathbf{H}\in \mathcal{V}_{n,K}$ and $\mathbf{D}:n\times n$,
$\mathbf{D}=\diag (D_{1},\ldots,D_{n})$, $n=\min (N-1),K$; $D_{1}\geq D_{2}\geq
\cdots \geq D_{n}\geq 0$, such that $\mathbf{Y}=\mathbf{V}'\mathbf{D}\mathbf{H}$;
This factorization is termed non-singular part of the SVD. Then
\begin{equation*}
    (d\mathbf{Y})=2^{-n}|\mathbf{D}|^{N-1+K-2n}\prod_{i<j}^{n}(D_{i}^{2}-D_{j}^{2})(d\mathbf{D})
    (\mathbf{V}d\mathbf{V}')(\mathbf{H}d\mathbf{H}').
\end{equation*}
\end{lem}
\textit{Proof.} See \citet{dgm:97}. \findemo

So, we can obtain:

\begin{thm}\label{th:jointVD}
The joint density of $(\mathbf{V},\mathbf{D})$ is
\begin{eqnarray*}
   f_{\mathbf{V},\mathbf{D}}(\mathbf{V},\mathbf{D})&=&\frac{\pi^{\frac{nK}{2}}|\mathbf{D}|^{N-1+K-2n}
   \prod_{i<j}(D_{i}^{2}-D_{j}^{2})}{\Gamma_{n}\left[\frac{K}{2}\right]|\boldsymbol{\Sigma}|^{\frac{K}{2}}}\\
   &&\times\sum_{t=0}^{\infty}\sum_{\kappa}\frac{h^{(2t)}\left[\tr\left(\boldsymbol{\Sigma}^{-1}
   \mathbf{V}'\mathbf{D}^{2}\mathbf{V}+\boldsymbol{\Omega}\right)\right]} {t!\left(\frac{1}{2}K\right)_{\kappa}}
   C_{\kappa}\left(\boldsymbol{\Omega}\boldsymbol{\Sigma}^{-1}\mathbf{V}'\mathbf{D}^{2}\mathbf{V}\right).
\end{eqnarray*}
\end{thm}
\textit{Proof.} Let be $\boldsymbol{\Omega} = \boldsymbol{\Sigma}^{-1}
\boldsymbol{\mu}\boldsymbol{\Theta}^{-1}\boldsymbol{\mu}'$, so the density of
$\mathbf{Y}$ is given by
$$
  f_{\mathbf{Y}}(\mathbf{Y})=\frac{1}{|\boldsymbol{\Sigma}|^{\frac{K}{2}}}h\left[\tr\left(\boldsymbol{\Sigma}^{-1}
  \mathbf{Y}\mathbf{Y}'+\boldsymbol{\Omega}\right)-2\tr\boldsymbol{\mu}'\boldsymbol{\Sigma}^{-1}\mathbf{Y}\right].
$$
Now, make the change of variables
$\mathbf{Y}=\mathbf{V}'\mathbf{D}\mathbf{H}$, so, by Lemma
\ref{lem:SVDjacobian}, the joint density function of $\mathbf{V}$,
$\mathbf{D}$, $\mathbf{H}$ is
\begin{eqnarray*}
    dF_{\mathbf{V},\mathbf{D},\mathbf{H}}(\mathbf{V},\mathbf{D},\mathbf{H})&=&
    \frac{2^{-n}|\mathbf{D}|^{N-1+K-2n}\displaystyle\prod_{i<j}(D_{i}^{2}-D_{j}^{2})}{|\boldsymbol{\Sigma}|^{\frac{K}{2}}}
    (\mathbf{V}d\mathbf{V}')(d\mathbf{D})\\
    && \times  h\left[\tr\left(\boldsymbol{\Sigma}^{-1}\mathbf{V}'\mathbf{D}^{2}\mathbf{V} +
    \boldsymbol{\Omega}\right)-2\tr\boldsymbol{\mu}'\boldsymbol{\Sigma}^{-1}\mathbf{V}'\mathbf{D}
    \mathbf{H}\right](\mathbf{H}d\mathbf{H}').
\end{eqnarray*}
Expanding in power series
\begin{eqnarray*}
    &&dF_{\mathbf{V},\mathbf{D},\mathbf{H}}(\mathbf{V},\mathbf{D},\mathbf{H})=\frac{2^{-n}
    |\mathbf{D}|^{N-1+K-2n}\displaystyle\prod_{i<j}(D_{i}^{2}-D_{j}^{2})}{|\boldsymbol{\Sigma}|^{\frac{K}{2}}}
    (\mathbf{V}d\mathbf{V}')(d\mathbf{D})\\
    && \quad\times \sum_{t=0}^{\infty}\frac{1}{t!} h^{(t)}\left[\tr\left(\boldsymbol{\Sigma}^{-1}
    \mathbf{V}'\mathbf{D}^{2}\mathbf{V}+\boldsymbol{\Omega}\right)\right]
    \left[\tr\left(-2\tr\boldsymbol{\mu}'\boldsymbol{\Sigma}^{-1}\mathbf{V}'\mathbf{D}\mathbf{H}\right)
    \right]^{t}(\mathbf{H}d\mathbf{H}').
\end{eqnarray*}
From Lemma \ref{lem:geneq22James64}
\begin{footnotesize}
\begin{equation*}
    \int_{\mathcal{V}_{n,K}}\left[\tr\left(-2\tr\boldsymbol{\mu}'\boldsymbol{\Sigma}^{-1}\mathbf{V}'
    \mathbf{D}\mathbf{H}\right)\right]^{2t}(\mathbf{H}d\mathbf{H}')=\frac{2^{n}\pi^{\frac{nK}{2}}}{\Gamma_{n}
    \left[\frac{K}{2}\right]} \sum_{\kappa}\frac{\left(\frac{1}{2}\right)_{t}4^{t}}{\left(\frac{1}{2}K\right)_{\kappa}}
    C_{\kappa}\left(\boldsymbol{\Omega}\boldsymbol{\Sigma}^{-1}\mathbf{V}'\mathbf{D}^{2}\mathbf{V}\right).
\end{equation*}
\end{footnotesize}
Observing that
$\displaystyle\frac{\left(\frac{1}{2}\right)_{t}4^{t}}{(2t)!}=\frac{1}{t!}$, the
marginal joint density of $\mathbf{V}$, $\mathbf{D}$ is given by
\begin{eqnarray*}
    &&dF_{\mathbf{V},\mathbf{D}}(\mathbf{V},\mathbf{D})=\frac{\pi^{\frac{nK}{2}}|\mathbf{D}|^{N-1+K-2n}
    \displaystyle\prod_{i<j}(D_{i}^{2}-D_{j}^{2})} {\Gamma_{n}\left(\frac{K}{2}\right)
    |\boldsymbol{\Sigma}|^{\frac{K}{2}}} \\
    && \quad\times \sum_{t=0}^{\infty}\frac{h^{(2t)}\left[\tr\left(\boldsymbol{\Sigma}^{-1}\mathbf{V}'
    \mathbf{D}^{2}\mathbf{V}+\boldsymbol{\Omega}\right)\right]}{t!\left(\frac{1}{2}K\right)_{\kappa}}
    C_{\kappa}\left(\boldsymbol{\Omega}\boldsymbol{\Sigma}^{-1}\mathbf{V}'\mathbf{D}^{2}
    \mathbf{V}\right)(\mathbf{V}d\mathbf{V}')(d\mathbf{D}). \qed
\end{eqnarray*}
Now, let be
\begin{equation}\label{eq:R}
    \mathbf{R}'=\mathbf{V}'\mathbf{D}
\end{equation}
then $\mathbf{V}'=\mathbf{R}'\mathbf{D}^{-1}$, so $d\mathbf{V}=d\mathbf{R}'
\mathbf{D}^{-1}$ and
$\mathbf{V}d\mathbf{V}'=\mathbf{D}^{-1}\mathbf{R}d\mathbf{R}'\mathbf{D}^{-1}$. But
$\mathbf{R}d\mathbf{R}':n\times n$ is skew symmetric, thus
\begin{equation}\label{eq:jacobianR}
    (\mathbf{V}d\mathbf{V}')=|\mathbf{D}|^{-n+1}(\mathbf{R}d\mathbf{R}')
\end{equation}

\begin{thm}\label{th:SVDreflectionsizeandshape}
The SVD reflection size-and-shape density is
\begin{eqnarray}\label{eq:SVDreflectionsizeandshape}
    dF_{\mathbf{R}}(\mathbf{R})&=&\frac{2^{-n}\Gamma_{n}\left[\frac{N+K-2n}{2}\right]\pi^{\frac{nK}{2}}}
    {\pi^{n(N+K-(3n-1)/2)/2}\Gamma_{n}\left[\frac{K}{2}\right]|\boldsymbol{\Sigma}|^{\frac{K}{2}}}(d\mathbf{D})\nonumber\\
    && \times \sum_{t=0}^{\infty}\frac{h^{(2t)}\left[\tr\left(\boldsymbol{\Sigma}^{-1}\mathbf{R}'\mathbf{R}+
    \boldsymbol{\Omega}\right)\right] C_{\kappa}\left(\boldsymbol{\Omega}\boldsymbol{\Sigma}^{-1}\mathbf{R}'
    \mathbf{R}\right)}{t!\left(\frac{1}{2}K\right)_{\kappa}}(\mathbf{R}d\mathbf{R}').
\end{eqnarray}
\end{thm}
\textit{Proof.}  The joint density function of $\mathbf{V},
\mathbf{D}$ is
\begin{small}
\begin{eqnarray*}
    dF_{\mathbf{V},\mathbf{D}}(\mathbf{V},\mathbf{D})&=&\frac{\pi^{\frac{nK}{2}}|\mathbf{D}|^{N-1+K-2n}
    \displaystyle\prod_{i<j}(D_{i}^{2}-D_{j}^{2})}{\Gamma_{n}\left[\frac{K}{2}\right]
    |\boldsymbol{\Sigma}|^{\frac{K}{2}}}(d\mathbf{D}) \\
    && \times\sum_{t=0}^{\infty}\sum_{\kappa}\frac{h^{(2t)}\left[\tr\left(\boldsymbol{\Sigma}^{-1}
    \mathbf{V}'\mathbf{D}^{2}\mathbf{V}+\boldsymbol{\Omega}\right)\right]
    C_{\kappa}\left(\boldsymbol{\Omega}\boldsymbol{\Sigma}^{-1}\mathbf{V}'
    \mathbf{D}\mathbf{V}\right)}{t!\left(\frac{1}{2}K\right)_{\kappa}}\left(\mathbf{V}d\mathbf{V}'\right).
\end{eqnarray*}
\end{small}
Now, let be $\mathbf{R}'=\mathbf{V}'\mathbf{D}$ then
$\mathbf{V}'=\mathbf{R}'\mathbf{D}^{-1}$, so
$d\mathbf{V}'=d\mathbf{R}'\mathbf{D}^{-1}$ and
$\mathbf{V}d\mathbf{V}'=\mathbf{D}^{-1}\mathbf{R}d\mathbf{R}'\mathbf{D}^{-1}$, where
$\mathbf{R}d\mathbf{R}'$ is an $n\times n$ skew-symmetric matrix.
$$
  \left(\mathbf{V}d\mathbf{V}'\right)=|\mathbf{D}^{-1}|^{n-1}\left(\mathbf{R}d\mathbf{R}'
  \right)=|\mathbf{D}|^{-n+1}\left(\mathbf{R}d\mathbf{R}'\right).
$$
Thus the joint density function of $\mathbf{R}$, $\mathbf{D}$ is
\begin{eqnarray*}
    dF_{\mathbf{R},\mathbf{D}}(\mathbf{R},\mathbf{D})&=&\frac{\pi^{\frac{nK}{2}}|\mathbf{D}|^{N+K-3n}
    \displaystyle\prod_{i<j}(D_{i}^{2}-D_{j}^{2})}{\Gamma_{n}\left[\frac{K}{2}\right]
    |\boldsymbol{\Sigma}|^{\frac{K}{2}}}(d\mathbf{D})\\
    && \times\sum_{t=0}^{\infty}\sum_{\kappa}\frac{h^{(2t)}\left[\tr\left(\boldsymbol{\Sigma}^{-1}
    \mathbf{R}'\mathbf{R}+\boldsymbol{\Omega}\right)\right] C_{\kappa}\left(\boldsymbol{\Omega}
    \boldsymbol{\Sigma}^{-1}\mathbf{R}'\mathbf{R}\right)}{t!\left(\frac{1}{2}K\right)_{\kappa}}
    \left(\mathbf{R}d\mathbf{R}'\right).
\end{eqnarray*}
For the integration with respect to $\mathbf{D}$ note that if
$\mathbf{D}^{2}=\mathbf{L}$, so $dL_{i}=2D_{i}dD_{i}$ and
$(dL)=2^{n}|\mathbf{D}|(d\mathbf{D})$, thus
$$
  (d\mathbf{D})=2^{-n}|\mathbf{L}|^{-\frac{1}{2}}(d\mathbf{L}).
$$
Let be
$$
  J=\int_{\mathbf{D}}|\mathbf{D}|^{N+K-3n}\prod_{i<j}\left(D_{i}^{2}-D_{j}^{2}\right)(d\mathbf{D})
$$
so
\begin{eqnarray*}
    J&=&\int_{\mathbf{L}}|\mathbf{L}^{\frac{1}{2}}|^{N+K-3n}\prod_{i<j}\left(L_{i}-L_{j}\right)2^{-n}
    |\mathbf{L}|^{-\frac{1}{2}}(d\mathbf{L})\\
    &=&2^{-n}\int_{\mathbf{L}}|\mathbf{L}|^{\frac{1}{2}(N+K-3n-1)}\prod_{i<j}\left(L_{i}-L_{j}\right)(d\mathbf{L})
\end{eqnarray*}
From \citet{fz:90}, eq. (3.29), p.102.  we have that
\begin{eqnarray*}
    J&=&\frac{\displaystyle\prod_{i=1}^{n}\Gamma\left[\frac{1}{2}\left(N+k-2n-i+1\right)\right]}{\pi^{\frac{1}{2}n(N+K-2n)
    + \frac{1}{2}n}} \frac{\pi^{\frac{1}{2}n(N+K-2n)+\frac{1}{2}}}{\displaystyle\prod_{i=1}^{n}\Gamma\left[\frac{1}{2}
    \left(N+K-2n-i+1\right)\right]}\\
    && \times\int_{\mathbf{L}}|\mathbf{L}|^{[(N+K-2n)-n-1]/2}\prod_{i<j}\left(L_{i}-L_{j}\right)(d\mathbf{L})\\
    &=&\frac{2^{-n}\pi^{n(n-1)/4}\displaystyle\prod_{i=1}\Gamma\left[\frac{1}{2}(N+K-2n-i+1)\right]}{\pi^{\frac{n}{2}(N+K-2n)
    \frac{n}{2}+\frac{n(n-1)}{4}}}\\
    &=&\frac{\Gamma_{n}\left[\frac{1}{2}(N+K-2n)\right]2^{-n}}{\pi^{\frac{n}{2}\left(N+K-\frac{3n-1}{2}\right)}},
\end{eqnarray*}
then
\begin{eqnarray*}
    dF_{\mathbf{R}}(\mathbf{R})&=&\frac{\Gamma_{n}\frac{N+K-2n}{2}2^{-n}\pi^{\frac{nK}{2}}}{\Gamma_{n}
    \left[\frac{K}{2}\right] \pi^{\frac{n}{2}\left(N+K-\frac{3n-1}{2}\right)}
    |\boldsymbol{\Sigma}|^{\frac{K}{2}}}\\
    && \times \sum_{t=0}^{\infty}\sum_{\kappa}\frac{h^{(2t)}\left[\tr\left(\boldsymbol{\Sigma}^{-1}
    \mathbf{R}'\mathbf{R}+\boldsymbol{\Omega}\right)\right]  C_{\kappa}\left(\boldsymbol{\Omega}
    \boldsymbol{\Sigma}^{-1}\mathbf{R}'\mathbf{R}\right)}{t!\left(\frac{1}{2}K\right)_{\kappa}}
    (\mathbf{R}d\mathbf{R}') \qed
\end{eqnarray*}

\section{Reflection Shape Density}\label{sec:svdreflectionshape}

For the SVD reflection shape density consider the following transformations
$$
  \mathbf{L}\mathbf{X} \boldsymbol{\Theta}^{-1/2}=
  \mathbf{L}\mathbf{Z}=\mathbf{Y}=(\mathbf{V}'\mathbf{D})\mathbf{H}\equiv
  \mathbf{R}'\mathbf{H}=r\mathbf{W}\mathbf{H}=r\mathbf{W}(\mathbf{u})\mathbf{H},
$$
where $\mathbf{R}'=\mathbf{V}'\mathbf{D}$ and $\mathbf{W} = \mathbf{R}'/r$.

Now, note that $\mathbf{V}'\mathbf{D}$ contains $(N-1)n$ coordinates. Then
$$
  \vec \mathbf{W}=\frac{1}{r}\vec (\mathbf{V}'\mathbf{D}), \quad
  r=\|\mathbf{V}'\mathbf{D}\|=\sqrt{\tr
  \mathbf{V}'\mathbf{D}^{2}\mathbf{V}}=\|\mathbf{Y}\|.
$$
Then, by \citet[Theorem 2.1.3, p. 55]{MR1982}:
\begin{eqnarray*}
    (d\vec \mathbf{W}(\mathbf{u}))&=&r^{m}\prod_{i=1}^{m}\sin^{m-i}\theta_{i}\bigwedge_{i=1}^{m}
    d\theta_{i}\wedge dr\\
    &=& r^{m}J(\mathbf{u})\bigwedge_{i=1}^{m}d\theta_{i}\wedge dr,
\end{eqnarray*}
with $m=(N-1)n-1$, $\mathbf{u}=(\theta_{1},\ldots,\theta_{m})'$.

Hence,
\begin{thm}\label{th:SVDreflectionshape}
The SVD reflection shape density is given by
\begin{eqnarray}\label{eq:SVDreflectionshape}
    dF_{\mathbf{W}}(\mathbf{W})&=&\frac{\Gamma_{n}\left[\frac{N+K-2n}{2}\right]J(\mathbf{u})
    \pi^{\frac{nK}{2}}}{2^{n}\pi^{\frac{n}{2}\left(N+K-\frac{3n-1}{2}\right)}
    \Gamma_{n}\left[\frac{K}{2}\right]|\boldsymbol{\Sigma}|^{\frac{K}{2}}}
    \sum_{t=0}^{\infty}\sum_{\kappa}\frac{C_{\kappa}\left(\boldsymbol{\Omega}
    \boldsymbol{\Sigma}^{-1}\mathbf{W}'\mathbf{W}\right)}{t!\left(\frac{1}{2}K\right)_{\kappa}}
    (\mathbf{W}d\mathbf{W}') \nonumber\\
    && \times\int_{0}^{\infty}r^{m+n+2t-1}h^{(2t)}\left[r^{2}\tr\boldsymbol{\Sigma}^{-1}
    \mathbf{W}'\mathbf{W}+\tr\boldsymbol{\Omega}\right](dr).
\end{eqnarray}
\end{thm}
\textit{Proof.} The density of $\mathbf{R}$ is
\begin{eqnarray*}
    dF_{\mathbf{R}}(\mathbf{R})&=&\frac{\Gamma_{n}\left[\frac{N+K-2n}{2}\right]2^{-n}
    \pi^{\frac{nK}{2}}}{\pi^{\frac{n}{2}\left[N+K-\frac{3n-1}{2}\right]
    }\Gamma_{n}\left[\frac{K}{2}\right]|\boldsymbol{\Sigma}|^{\frac{K}{2}}}\\
    && \times\sum_{t=0}^{\infty}\sum_{\kappa}\frac{h^{(2t)}\left[\tr\left(\boldsymbol{\Sigma}^{-1}
    \mathbf{R}'\mathbf{R}+\boldsymbol{\Omega}\right)\right]C_{\kappa}\left(\boldsymbol{\Omega}
    \boldsymbol{\Sigma}^{-1}\mathbf{R}'\mathbf{R}\right)}
    {t!\left(\frac{1}{2}K\right)_{\kappa}}(\mathbf{R}d\mathbf{R}').
\end{eqnarray*}
Putting $\mathbf{W}(\mathbf{u})=\mathbf{R}'/r$, the joint density of $r$ and
$\mathbf{u}$ is\\[2ex]
$
    dF_{r,\mathbf{W}(\mathbf{u})}(r,\mathbf{W}(\mathbf{u}))  = \displaystyle\frac{\Gamma_{n}
    \left[\frac{N+K-2n}{2}\right]2^{-n}
    \pi^{\frac{nK}{2}}r^{m}J(\mathbf{u})} {\pi^{\frac{n}{2}\left[N+K-\frac{3n-1}{2}\right]
    }\Gamma_{n}\left[\frac{K}{2}\right]|\boldsymbol{\Sigma}|^{\frac{K}{2}}}\left(r^{2}
    \mathbf{W}d\mathbf{W}'\right)
$
\par \noindent \hfill
    \hbox{$\displaystyle\times\sum_{t=0}^{\infty}\sum_{\kappa}
    \frac{h^{(2t)}\left[\tr\left(r^{2}\boldsymbol{\Sigma}^{-1}\mathbf{W}'\mathbf{W}+
    \boldsymbol{\Omega}\right)\right]C_{\kappa}\left(r^{2}\boldsymbol{\Omega}
    \boldsymbol{\Sigma}^{-1}\mathbf{W}'\mathbf{W}\right)}
    {t!\left(\frac{1}{2}K\right)_{\kappa}}.$}
\par \noindent

Note that
\begin{enumerate}
    \item
    $C_{\kappa}\left(r^{2}\boldsymbol{\Omega}\boldsymbol{\Sigma}^{-1}\mathbf{W}'\mathbf{W}\right)=
    r^{2t}C_{\kappa}\left(\boldsymbol{\Omega}\boldsymbol{\Sigma}^{-1}\mathbf{W}'\mathbf{W}\right)$.
    \item
    $\left(r^{2}\mathbf{W}d\mathbf{W}'\right)=\left(\left(r\mathbf{I}\right)\mathbf{W}d\mathbf{W}'
    \left(r\mathbf{I}\right)\right)=|r\mathbf{I}|^{n-1}\left(\mathbf{W}d\mathbf{W}'\right)=r^{n-1}
    \left(\mathbf{W}d\mathbf{W}'\right)$.
\end{enumerate}
Collecting powers of $r$, the marginal of $\mathbf{W}$ is
\begin{eqnarray*}
    dF_{\mathbf{W}}(\mathbf{W})&=&\frac{\Gamma_{n}\left[\frac{N+K-2n}{2}\right]2^{-n}\pi^{\frac{nK}{2}}J(\mathbf{u})}
    {\pi^{\frac{n}{2}\left[N+K-\frac{3n-1}{2}\right]
    }\Gamma_{n}\left[\frac{K}{2}\right]|\boldsymbol{\Sigma}|^{\frac{K}{2}}}\left(\mathbf{W}d\mathbf{W}'\right)
    \sum_{t=0}^{\infty}\sum_{\kappa}
    \frac{C_{\kappa}\left(\boldsymbol{\Omega}\boldsymbol{\Sigma}^{-1}\mathbf{W}'\mathbf{W}\right)}
    {t!\left(\frac{1}{2}K\right)_{\kappa}}\\&&\times\int_{0}^{\infty}r^{m+n+2t-1}h^{(2t)}
    \left[r^{2}\tr\boldsymbol{\Sigma}^{-1}\mathbf{W}'\mathbf{W}+\tr\boldsymbol{\Omega}\right](dr).
    \qed
\end{eqnarray*}

\section{Distributions excluding reflection}\label{sec:SVDexcludingreflection}

Recall that  the SVD shape coordinates ${\bf u}$ of $\mathbf{X}$ are obtained as
follows
\begin{equation}
    \mathbf{L}\mathbf{X} \boldsymbol{\Theta}^{-1/2}=
    \mathbf{L}\mathbf{Z} = \mathbf{Y} = \mathbf{H}'\mathbf{D}\mathbf{P}
    = r \mathbf{W}\mathbf{P} = r\mathbf{W}({\bf u})\mathbf{P}
\end{equation}
where the SVD shape coordinate system is given by $\mathbf{H}'\mathbf{D}$ (\citet
[pp. 296-298]{g:91}) and $r =||\mathbf{H}\mathbf{D}|| = (\tr \mathbf{D} \mathbf{H}
\mathbf{H}'\mathbf{D})^{1/2} = (\tr \mathbf{D}^{2})^{1/2} = ||\mathbf{D}|| =
||\mathbf{Y}||$. When $n = K $ we studied the distributions including reflection,
i.e.  $\mathbf{P} \in \mathcal{O}(k)$, $|\mathbf{P}| = \pm 1$, $\mathbf{D}_{K} \geq
0$ and $(\mathbf{H},\mathbf{D})$, written $(\mathbf{H},\mathbf{D})^{R}$ for
definiteness, contains reflection SVD shape co-ordinates.

In this section we consider the case when  $\mathbf{P}$ excludes reflection, thus
$\mathbf{P} \in \mathcal{SO}(K)$, $|\mathbf{P}| = + 1$, $|D_{K}| \geq 0$,
$\sign(D_{m}) = \sign|\mathbf{X}|$ and $(\mathbf{H},\mathbf{D})$, may be written
$(\mathbf{H},\mathbf{D})^{NR}$ for definiteness.

Finally, we have that the  excluding reflection SVD size-and-shape and SVD shape
densities are given by (\ref{eq:SVDreflectionsizeandshape}) and
(\ref{eq:SVDreflectionshape}) divided by 2,  respectively.

\section{Central Case}\label{sec:SVDcentralcase}

Now, we can derive easily the corresponding central distributions of this work.

\begin{cor}\label{cor:SVDCentralreflectionsizeandshape}
The central reflection SVD size-and-shape density is
$$
  dF_{\mathbf{R}}(\mathbf{R})=\frac{2^{-n}\Gamma_{n}\left[\frac{N+K-2n}{2}\right]\pi^{\frac{nK}{2}}}
  {\pi^{\frac{n}{2}\left[N+K-\frac{3n-1}{2}\right]}\Gamma_{n}\left[\frac{K}{2}\right]
  |\boldsymbol{\Sigma}|^{\frac{K}{2}}}  h\left[\tr\boldsymbol{\Sigma}^{-1}\mathbf{R}'\mathbf{R}\right]
  (\mathbf{R}d\mathbf{R}').
$$
\end{cor}
\textit{Proof.} Just take $\boldsymbol{\mu}=0$ in Theorem
\ref{th:SVDreflectionsizeandshape} and recall that $h^{(0)}(\cdot)\equiv h(\cdot)$.
\findemo

Finally,
\begin{cor}\label{cor:SVDCentralreflectionshapeINVARIANT}
The central reflection SVD shape density is invariant under the elliptical family and
it is given by
$$
  dF_{\mathbf{W}}(\mathbf{W})=\frac{2^{-n-1}\Gamma_{n}\left[\frac{N+K-2n}{2}\right]\pi^{\frac{nK}{2}}
  \Gamma\left[\frac{m+n}{2}\right]}{\pi^{\frac{n}{2}\left[N+K-\frac{3n-1}{2}\right] + \frac{m+n}{2}}
  \Gamma_{n}\left[\frac{K}{2}\right]|\boldsymbol{\Sigma}|^{\frac{K}{2}}}
  h\left[\tr\boldsymbol{\Sigma}^{-1}\mathbf{W}'\mathbf{W}\right](\mathbf{W}d\mathbf{W}').
$$
\end{cor}
\textit{Proof.} Taking $\boldsymbol{\mu}=0$ and $s=\left(\tr
\boldsymbol{\Sigma}^{-1}\mathbf{W}'\mathbf{W}\right)^{1/2}r$ in Theorem
\ref{th:SVDreflectionshape} we obtain the result, since\\[2ex]
$
  \displaystyle\int_{0}^{\infty}r^{m+n-1}h\left[r^{2}\tr\boldsymbol{\Sigma}^{-1}\mathbf{W}'\mathbf{W}\right]dr
$
\begin{eqnarray*}
    \phantom{XXXXXXXXXX}
    &=& \int_{0}^{\infty}\left(\frac{s}{\left(\tr\boldsymbol{\Sigma}^{-1}\mathbf{W}'
    \mathbf{W}\right)^{\frac{1}{2}}}\right) h\left(s^{2}\right)\frac{ds}{\left(\tr\boldsymbol{\Sigma}^{-1}
    \mathbf{W}'\mathbf{W}\right)^{\frac{1}{2}}}\\
    &=& \left(\tr\boldsymbol{\Sigma}^{-1}\mathbf{W}'\mathbf{W}\right)^{-\frac{m+n}{2}}\frac{\Gamma
    \left[\displaystyle\frac{m+n}{2}\right]}{2\pi^{\frac{m+n}{2}}}.
\end{eqnarray*}

\section{Some particular models}\label{sec:particularmodels}

Finally, we give explicit shapes densities for some elliptical models.

The Kotz type I model is given by
$$
  h(y) = \frac{R^{T-1+\frac{K(N-1)}{2}}\Gamma\left[\frac{K(N-1)}{2}\right]}{\pi^{K(N-1)/2}
  \Gamma\left[T-1+\frac{K(N-1)}{2}\right]}y^{T-1}\exp\{-Ry\}.
$$
Then, the corresponding $k$-th derivative of $h$, follows from
$\displaystyle\frac{d^{k}}{dy^{k}}y^{T-1}\exp\{-Ry\}$, which is given by
$$
  (-R)^{k}y^{T-1}\exp\{-Ry\}\left\{1+\sum_{m=1}^{k}\binom{k}{m}
  \left[\prod_{i=0}^{m-1}(T-1-i)\right](-Ry)^{-m}\right\},\nonumber\\
$$
see \citet{Caro2009}.

It is of interest the Gaussian case, i.e. when $T=1$ and $R=\frac{1}{2}$, here the
derivation is straightforward from the general density.

The required derivative follows easily, it is,
$h^{(k)}(y)=\displaystyle\frac{R^{\frac{K(N-1)}{2}}}{\pi^{\frac{K(N-1)}{2}}}(-R)^{k}\exp\{-Ry\}$
and
\begin{eqnarray*}
    \int_{0}^{\infty}
    &&r^{m+n+2t-1}h^{(2t)}\left[r^{2}\tr\boldsymbol{\Sigma}^{-1}\mathbf{W}'\mathbf{W}+
    \tr\boldsymbol{\Omega}\right](dr) \\
    && = \frac{R^{\frac{M}{2}-\frac{1}{2}(m+n)+t}}{2\pi^{\frac{M}{2}}} \exp\{-R \tr\boldsymbol{\Omega}\}\left(\tr
    \boldsymbol{\Sigma}^{-1}\mathbf{W}\mathbf{W}'\right)^{-\frac{m+n}{2}-t}\Gamma\left[\frac{m+n}{2}+t\right].
\end{eqnarray*}
Hence $dF_{\mathbf{W}}(\mathbf{W})$ is given by
\begin{eqnarray*}
    &=&\frac{\Gamma_{n}\left[\frac{N+K-2n}{2}\right]J(\mathbf{u})
    \pi^{\frac{nK}{2}}}{2^{n}\pi^{\frac{n}{2}\left(N+K-\frac{3n-1}{2}\right)}
    \Gamma_{n}\left[\frac{K}{2}\right]|\boldsymbol{\Sigma}|^{\frac{K}{2}}}
    \sum_{t=0}^{\infty}\sum_{\kappa}\frac{C_{\kappa}\left(\boldsymbol{\Omega}
    \boldsymbol{\Sigma}^{-1}\mathbf{W}'\mathbf{W}\right)}{t!\left(\frac{1}{2}K\right)_{\kappa}}
    (\mathbf{W}d\mathbf{W}')\\
    && \times\int_{0}^{\infty}r^{m+n+2t-1}h^{(2t)}\left[r^{2}\tr\boldsymbol{\Sigma}^{-1}\mathbf{W}'
    \mathbf{W}+\tr\boldsymbol{\Omega}\right](dr)\\
    &=&\frac{R^{\frac{M}{2}-\frac{1}{2}(m+n)}\Gamma_{n}\left[\frac{N+K-2n}{2}\right]J(\mathbf{u})\exp\{-R
    \tr\boldsymbol{\Omega}\}}{2^{n+1}\pi^{\frac{n}{2}\left(N-\frac{3n-1}{2}\right)+\frac{M}{2}}
    \Gamma_{n}\left[\frac{K}{2}\right]|\boldsymbol{\Sigma}|^{\frac{K}{2}}}
    \sum_{t=0}^{\infty}\sum_{\kappa}\frac{C_{\kappa}\left(R\boldsymbol{\Omega}\boldsymbol{\Sigma}^{-1}
    \mathbf{W}'\mathbf{W}\right)}{t!\left(\frac{1}{2}K\right)_{\kappa}} \\
    && \times \left(\tr \boldsymbol{\Sigma}^{-1}\mathbf{W}\mathbf{W}'\right)^{-\frac{m+n}{2}-t}
    \Gamma\left[\frac{m+n}{2}+t\right](\mathbf{W}d\mathbf{W}').
\end{eqnarray*}
Therefore, we have proved that

\begin{cor}\label{th:SVDreflectionshapeNORMAL}
The Gaussian SVD reflection shape density is
\begin{small}
\begin{eqnarray}\label{eq:SVDreflectionshapeNormal}
\hspace{-0.75cm}
    dF_{\mathbf{W}}(\mathbf{W})&=&\frac{R^{\frac{M}{2}-\frac{1}{2}(m+n)}\Gamma_{n}\left[\frac{N+K-2n}{2}\right]
    J(\mathbf{u})\etr\{-R\boldsymbol{\Omega}\}}{2^{n+1}\pi^{\frac{n}{2}\left(N-\frac{3n-1}{2}\right)+\frac{M}{2}}
    \Gamma_{n}\left[\frac{K}{2}\right]|\boldsymbol{\Sigma}|^{\frac{K}{2}}}
    \sum_{t=0}^{\infty}\sum_{\kappa}\frac{C_{\kappa}\left(R\boldsymbol{\Omega}\boldsymbol{\Sigma}^{-1}
    \mathbf{W}'\mathbf{W}\right)}{t!\left(\frac{1}{2}K\right)_{\kappa}}
    \nonumber\\
    &&\times \left(\tr \boldsymbol{\Sigma}^{-1}\mathbf{W}\mathbf{W}'\right)^{-\frac{m+n}{2}-t}
    \Gamma\left[\frac{m+n}{2}+t\right](\mathbf{W}d\mathbf{W}'),
\end{eqnarray}
\end{small}
where $M=(N-1)K$.
\end{cor}

Finally, we propose the result for the Kotz type I model
$$
  h(y) = \frac{R^{T-1+\frac{K(N-1)}{2}}\Gamma\left[\frac{K(N-1)}{2}\right]}{\pi^{K(N-1)/2}
  \Gamma\left[T-1+\frac{K(N-1)}{2}\right]}y^{T-1}\exp\{-Ry\},
$$

\begin{cor}\label{th:SVDreflectionshapeKotz}
The Kotz type I SVD reflection shape density is
\begin{eqnarray*}
    f_{\mathbf{W}}(\mathbf{W})&=&\frac{R^{T-1+\frac{M}{2}-\frac{m+n}{2}}\Gamma\left[\frac{M}{2}\right]
    \Gamma_{n}\left[\frac{N+K-2n}{2}\right]J(\mathbf{u})(\tr\boldsymbol{\Omega})^{T-1}
    \etr(-R\boldsymbol{\Omega})}{2^{n+1}\pi^{\frac{n}{2}\left(N-\frac{3n-1}{2}\right)+\frac{M}{2}}
    \Gamma_{n}\left[\frac{K}{2}\right]|\boldsymbol{\Sigma}|^{\frac{K}{2}}\Gamma\left[T-1+\frac{M}{2}\right]}\\
    && \times \sum_{t=0}^{\infty}\sum_{\kappa}\frac{C_{\kappa}\left(R\boldsymbol{\Omega}\boldsymbol{\Sigma}^{-1}
    \mathbf{W}'\mathbf{W}\right)}{t!\left(\frac{1}{2}K\right)_{\kappa}}(\mathbf{W}d\mathbf{W}')
    \left(\tr\boldsymbol{\Sigma}^{-1}\mathbf{W}'\mathbf{W}\right)^{-\frac{m+n}{2}-t}\\
    && \times\left\{\sum_{u=0}^{\infty}\frac{\Gamma\left[\frac{M}{2}+t+u\right]\displaystyle\prod_{s=0}^{u-1}(T-1-s)}
    {u!R^{u}(\tr\boldsymbol{\Omega})^{u} \Gamma\left[T-1+\frac{M}{2}\right] }\right. \\
    &&\left.+ \sum_{m=1}^{k}\binom{k}{m}\left[\prod_{i=0}^{m-1}(T-1-i)\right] \frac{(-R)^{-m}
    \left(\tr \boldsymbol{\Omega}\right)^{-m}}{\Gamma\left[T-1-m+\frac{M}{2}\right]}\right. \\
    && \times\left.\sum_{u=0}^{\infty}\frac{\Gamma\left[\frac{M}{2}+t+u\right]\displaystyle\prod_{s=0}^{u-1}
    (T-1-m-s)}{u!R^{u}(\tr\boldsymbol{\Omega})^{u}}\right\},
\end{eqnarray*}
where $M=(N-1)K$.
\end{cor}

\textit{Proof.} As we note before  the $k$-th derivative of $h$ follows from,
\begin{eqnarray*}
    &&\frac{d^{k}}{dy^{k}}y^{T-1}\exp\{-Ry\}\\&=&(-R)^{k}y^{T-1}\exp\{-Ry\}\left\{1+\sum_{m=1}^{k}
  \binom{k}{m}\left[\prod_{i=0}^{m-1}(T-1-i)\right](-Ry)^{-m}\right\},
\end{eqnarray*}
and the corresponding SVD reflection shape density, $dF_{\mathbf{W}}(\mathbf{W})$,
is obtained after some simplification as
\begin{eqnarray*}
    &=&\frac{\Gamma_{n}\left[\frac{N+K-2n}{2}\right]J(\mathbf{u})
    \pi^{\frac{nK}{2}}}{2^{n}\pi^{\frac{n}{2}\left(N+K-\frac{3n-1}{2}\right)}
    \Gamma_{n}\left[\frac{K}{2}\right]|\boldsymbol{\Sigma}|^{\frac{K}{2}}}
    \sum_{t=0}^{\infty}\sum_{\kappa}\frac{C_{\kappa}\left(\boldsymbol{\Omega}\boldsymbol{\Sigma}^{-1}
    \mathbf{W}'\mathbf{W}\right)}{t!\left(\frac{1}{2}K\right)_{\kappa}}(\mathbf{W}d\mathbf{W}')\\
    && \times\int_{0}^{\infty}r^{m+n+2t-1}h^{(2t)}\left[r^{2}\tr\boldsymbol{\Sigma}^{-1}\mathbf{W}'
    \mathbf{W}+\tr\boldsymbol{\Omega}\right](dr)\\
    &=&\frac{R^{T-1+\frac{M}{2}-\frac{m+n}{2}}\Gamma\left[\frac{M}{2}\right]\Gamma_{n}\left[\frac{N+K-2n}{2}
    \right]J(\mathbf{u})(\tr\boldsymbol{\Omega})^{T-1}\etr(-R\boldsymbol{\Omega})}{2^{n+1}\pi^{\frac{n}{2}
    \left(N-\frac{3n-1}{2}\right)+\frac{M}{2}}\Gamma_{n}\left[\frac{K}{2}\right]
    |\boldsymbol{\Sigma}|^{\frac{K}{2}}\Gamma\left[T-1+\frac{M}{2}\right]}\\
    &&\times \sum_{t=0}^{\infty}\sum_{\kappa}\frac{C_{\kappa}\left(R\boldsymbol{\Omega}
    \boldsymbol{\Sigma}^{-1}\mathbf{W}'\mathbf{W}\right)}{t!\left(\frac{1}{2}K\right)_{\kappa}}
    (\mathbf{W}d\mathbf{W}') \left(\tr\boldsymbol{\Sigma}^{-1}\mathbf{W}'\mathbf{W}\right)^{-\frac{m+n}{2}-t}\\
    && \times\left\{\sum_{u=0}^{\infty}\frac{\Gamma\left[\frac{M}{2}+t+u\right]\displaystyle\prod_{s=0}^{u-1}
    (T-1-s)}{u!R^{u}(\tr\boldsymbol{\Omega})^{u} \Gamma\left[T-1+\frac{M}{2}\right]}\right.\\
    &&\left.+ \sum_{m=1}^{k}\binom{k}{m}\left[\prod_{i=0}^{m-1}(T-1-i)\right] \frac{(-R)^{-m}\left(\tr
    \boldsymbol{\Omega}\right)^{-m}} {\Gamma\left[T-1-m+\frac{M}{2}\right]}\right.\\
    &&\times\left.\sum_{u=0}^{\infty}\frac{\Gamma\left[\frac{M}{2}+t+u\right]\displaystyle\prod_{s=0}^{u-1}
    (T-1-m-s)}{u!R^{u}(\tr\boldsymbol{\Omega})^{u}}\right\}. \qed
\end{eqnarray*}
The Gaussian case can be derived again  by taking $T=1$ in the above result.

\subsection{Example: Mouse Vertebra}\label{sub:mouse}

This classical application is studied in the Gaussian case by \citet{DM98}. Here we
consider again the same model and contrasted it, via the modified BIC criterion, with
two non Gaussian models.

The isotropic Gaussian shape density is given by

\begin{eqnarray}\label{eq:SVDreflectionshapeNORMAL}
    dF_{\mathbf{W}}(\mathbf{W})&=&\frac{2^{-\frac{1}{2}(2-m+M+n)}\Gamma_{n}\left[\frac{N+K-2n}{2}\right]
    J(\mathbf{u})}{\pi^{\frac{1}{4}(2M+n-3n^{2}+2nM)}\sigma^{-(m-M+n)}\Gamma_{n}\left[\frac{K}{2}\right]}
    \etr\left(-\frac{\boldsymbol{\mu}'\boldsymbol{\mu}}{2\sigma^{2}}\right)\nonumber\\
    &&\times\sum_{t=0}^{\infty}\frac{\Gamma\left[\frac{m+n}{2}+t\right]}{t!}\sum_{\kappa}
    \frac{C_{\kappa}\left(\frac{1}{2\sigma^{2}}\boldsymbol{\mu}'\mathbf{W}\mathbf{W}'
    \boldsymbol{\mu}\right)}{\left(\frac{1}{2}K\right)_{\kappa}}(\mathbf{W}d\mathbf{W}'),
\end{eqnarray}
where $M=K(N-1)$, $n=\min\{ (N-1),K\}$ and $m=(N-1)n-1$. Here we study three models,
the Gaussian shape (N), and the Kotz (K) model for $T=2$ and $T=3$.

The shape density associated to the Kotz model indexed by $T=2$, $R=\frac{1}{2}$ (and
$s=1$) is given by:

\begin{eqnarray}\label{eq:SVDreflectionshapeNORMAL2}
    dF_{\mathbf{W}}(\mathbf{W})&=&\frac{2^{-\frac{1}{2}(-m+M+n)}\Gamma_{n}\left[\frac{N+K-2n}{2}\right]
    J(\mathbf{u})}{\pi^{\frac{1}{4}(2M+n-3n^{2}+2nM)}\sigma^{-(m-M+n)}M\Gamma_{n}\left[\frac{K}{2}\right]}
    \etr\left(-\frac{\boldsymbol{\mu}'\boldsymbol{\mu}}{2\sigma^{2}}\right) \nonumber\\
    &&\times\sum_{t=0}^{\infty}\frac{\left(\tr\left(\frac{\boldsymbol{\mu}'\boldsymbol{\mu}}{2\sigma^{2}}
    \right)-2t\right)\Gamma\left[\frac{m+n}{2}+t\right] + \Gamma\left[\frac{m+n}{2}+t+1\right]}{t!}\nonumber\\
    &&\times\sum_{\kappa}\frac{C_{\kappa}\left(\frac{1}{2\sigma^{2}}\boldsymbol{\mu}'\mathbf{W}
    \mathbf{W}'\boldsymbol{\mu}\right)}{\left(\frac{1}{2}K\right)_{\kappa}}
    (\mathbf{W}d\mathbf{W}').
\end{eqnarray}
And the corresponding density, $dF_{\mathbf{W}}(\mathbf{W})$, for the Kotz model
$T=3$, is obtained as:
\begin{eqnarray*}
    &&=\frac{2^{-\frac{1}{2}(-2-m+M+n)}\Gamma_{n}\left[\frac{N+K-2n}{2}\right]
    J(\mathbf{u})}{\pi^{\frac{1}{4}(2M+n-3n^{2}+2nM)} \sigma^{-(m-M+n)}M(M+2)\Gamma_{n}
    \left[\frac{K}{2}\right]}\etr\left(-\frac{\boldsymbol{\mu}'\boldsymbol{\mu}}{2\sigma^{2}}\right)\\
    &&\times\sum_{t=0}^{\infty}\left\{\left[4t^{2}-2t-4t\tr\left(\frac{\boldsymbol{\mu}'
    \boldsymbol{\mu}}{2\sigma^{2}}\right) +\tr^{2}\left(\frac{\boldsymbol{\mu}'\boldsymbol{\mu}}{2\sigma^{2}}
    \right)\right]\Gamma\left[\frac{m + n}{2} + t\right] \right.\\
    && \left.+\left[-4t+2\tr\left(\frac{\boldsymbol{\mu}'\boldsymbol{\mu}}{2\sigma^{2}}\right)\right]
    \Gamma\left[\frac{m+n}{2}+t+1\right]+\Gamma\left[\frac{m+n}{2}+t+2\right]\right\}\\
    &&\times\sum_{\kappa} \frac{C_{\kappa}\left(\frac{1}{2\sigma^{2}}\boldsymbol{\mu}'\mathbf{W}
    \mathbf{W}'\boldsymbol{\mu}\right)}{t!\left(\frac{1}{2}K\right)_{\kappa}}(\mathbf{W}d\mathbf{W}').
\end{eqnarray*}
In order to decide which the elliptical model is the best one, different criteria
have been employed for the model selection. We shall consider a modification of the
BIC statistic as discussed in \citet{YY07}, and which was first achieved by
\citet{ri:78} in a coding theory framework. The modified BIC is given by:
$$
  BIC^{*}=-2\mathfrak{L}(\widetilde{\boldsymbol{\boldsymbol{\mu}}},\widetilde{\sigma}^{2},h)
  +n_{p}(\log(n+2)-\log 24),
$$
where
$\mathfrak{L}(\widetilde{\boldsymbol{\boldsymbol{\mu}}},\widetilde{\sigma}^{2},h)$ is
the maximum of the log-likelihood function, $n$ is the sample size and $n_{p}$ is the
number of parameters to be estimated for each particular shape density.

As proposed by \citet{kr:95} and \citet{r:95}, the following selection criteria have
been employed for the model selection.

\begin{table}[ht]  \centering \caption{Grades of evidence
corresponding to values of the $BIC^{*}$ difference.}\label{table2}
\medskip
\renewcommand{\arraystretch}{1}
\begin{center}
  \begin{tabular}{cl}
    \hline
    $BIC^{*}$ difference & Evidence\\
    \hline
    0--2 & Weak\\
    2--6 & Positive \\
    6--10 & Strong\\
    $>$ 10 & Very strong\\
    \hline
  \end{tabular}
\end{center}
\end{table}

Fixing the variance of the process as 50 (the maximum median variances of the two
samples), the maximum likelihood estimators for location parameters associated with
the small and large groups are summarized in the following table:

\medskip

\begin{table}[ht!]
  \centering
  \caption{The maximum likelihood estimators}\label{tb:1}
  \medskip
\begin{center}
\begin{small}\label{Tab4MR}
\begin{tabular}{|c|c|c|c|c|c|c|}
  \hline
  Group& $BIC^{*}$ & $\widetilde{\mu}_{11}$ &$\widetilde{\mu}_{12}$& $\widetilde{\mu}_{21}$ & $\widetilde{\mu}_{22}$
   & $\widetilde{\mu}_{31}$   \\
   & $\build{K:T=3}{G}{K:T=2}$&&&&&\\
  \hline
  Small & $\build{-39.6272
}{-5.9146}{-23.0250}$ & $\build{2.1250}{-1.9214}{-6.1682}$ &
$\build{-47.8016 }{-42.5338}{-45.0331}$
   & $\build{16.3513}{14.1761}{14.6035}$
  &$\build{ -3.9691}{-4.8190}{-6.4983}$    &$\build{ 26.9744}{24.0766}{25.5710}$
  \\ \hline
  Large& $\build{-9.2156}{24.5000
}{7.3880}$ &$\build{-34.4190}{-23.1834}{-42.9562}$
&$\build{-29.3792}{-32.9246}{-1.6455}$
&$\build{5.9227}{8.1455}{-3.6885}$
  & $\build{-13.7687}{-10.5612}{-13.7522}$   & $\build{20.9197}{22.5117}{3.6578}$   \\
  \hline
\end{tabular}
\end{small}
\end{center}
\end{table}

\begin{table}[ht!]
\begin{center}
\begin{footnotesize}
\begin{tabular}{|c|c|c|c|c|}
  \hline
   $\widetilde{\mu}_{32}$ &$\widetilde{\mu}_{41}$ & $\widetilde{\mu}_{42}$ &
   $\widetilde{\mu}_{51}$ & $\widetilde{\mu}_{52}$
    \\
  \hline
  $\build{ 2.1288 }{-0.2605}{-2.6121}$& $\build{4.9123}{4.8195}{5.5892}$ &$\build{ 5.8552}{4.7981}{4.6337}$
  & $\build{ -33.0646}{-29.2691}{-30.8241}$ & $\build{ 0.4382}{-0.7072}{6.0498}$
  \\ \hline
  $\build{-20.2846}{-12.8674}{-27.4382}$&$\build{6.1263}{5.0300}{5.2368}$ &$\build{1.1366}{2.3324}{-2.7585}$
  & $\build{-23.1221}{-26.0251}{-0.9942}$ & $\build{27.5996}{18.6865}{34.1886}$ \\
  \hline
\end{tabular}
\end{footnotesize}
\end{center}
\end{table}

According to the modified BIC criterion, the Kotz model with parameters $T=3$,
$R=\frac{1}{2}$ and $s=1$ is the most appropriate  among the three elliptical
densities for modeling the data. There is a very strong difference between the non
Gaussian and  the classical Gaussian model in this experiment.

Let $\boldsymbol{\boldsymbol{\mu}}_{1}$ and $\boldsymbol{\boldsymbol{\mu}}_{2}$ be
the mean shape of the small and large groups, respectively. We test equal mean shape
under the best model, and the likelihood ratio (based on
$-2\log\Lambda\approx\chi_{10}^{2}$) for the test
$H_{0}:\boldsymbol{\boldsymbol{\mu}}_{1}=\boldsymbol{\boldsymbol{\mu}}_{2}$ vs
$H_{a}:\boldsymbol{\boldsymbol{\mu}}_{1}\neq\boldsymbol{\boldsymbol{\mu}}_{2}$,
provides the p-value $0.84$, which means that there are extremely evidence that the
mean shapes of the two groups are equal if the variance of the experiment is fixed in
50 (the maximum median of the variances of the two samples), a deeper study of this
case is suggested, because the variance estimation was problematic in the performed
inference procedure for these data. We  highlight that our intention is to illustrate
the technique and performed inference with an exact likelihood efficiently computable
after modification of the algorithms given for hypergeometric series (\citet{KE06}).

A final comment, for any elliptical model we can obtain the SVD reflection model,
however a nontrivial problem appears, the $2t$-th derivative of the generator model,
which can be seen as a partition theory problem. For The general case of a Kotz model
($s\neq1$), and another models like Pearson II and VII, Bessel, Jensen-logistic, we
can use formulae for these derivatives given by \citet{Caro2009}. The resulting
densities have again a form of a generalized series of zonal polynomials which can be
computed efficiently after some modification of existing works for hypergeometric
series (see \citet{KE06}), thus the inference over an exact density can be performed,
avoiding the use of any asymptotic distribution, and the initial transformation
avoids the invariant polynomials of \citet{d:80}, and it lets the inclusion of any
correlation among landmarks.

\section*{Acknowledgment}

This research work was supported  by University of Medellin (Medellin, Colombia) and
Universidad Aut\'onoma Agraria Antonio Narro (M\'{e}xico),  joint grant No. 469,
SUMMA group. Also, the first author was partially supported  by CONACYT - M\'exico,
research grant no. \ 138713 and IDI-Spain, Grants No. \ FQM2006-2271 and
MTM2008-05785 and the paper was written during J. A. D\'{\i}az- Garc\'{\i}a's stay as
a visiting professor at the Department of Statistics and O. R. of the University of
Granada, Spain.

\end{document}